\documentclass[12pt,a4paper]{article}
\usepackage[german, french, english]{babel}
\usepackage[a4paper]{geometry}
\usepackage{amsmath,amssymb,amsthm, xspace}
\usepackage[matrix,arrow]{xy}

\title{The Cox Ring of a Del Pezzo Surface Has Rational Singularities}
\author{Oleg N. Popov\thanks{119992, GSP-2, Moscow, Leninskie Gory, Moscow
State University, Faculty of Mechanics and Mathematics, Department of
Higher Algebra. E-mail: \tt popov@mccme.ru.} }
\date{}

\newcommand{\OO}{\mathcal{O}}
\newcommand{\A}{\mathbb{A}}
\newcommand{\ZZ}{\mathbb{Z}}
\newcommand{\PP}{\mathbb{P}}
\newcommand{\GG}{\mathbb{G}}
\newcommand{\myT}{\mathfrak{T}}
\newcommand{\T}{\widetilde}
\newcommand{\WH}{\widehat}
\newcommand{\comment}[1]{}
\newcommand{\G}{\Gamma}
\newcommand{\wrt}{w.r.t.\@\xspace}
\newcommand{\ie}{i.e.,\xspace}
\newcommand{\co}{\colon}
\newcommand{\pp}{$/\!/$ }
\newcommand{\ang}{\mathring}

\newtheorem{stat}{Proposition}
\newtheorem{lemma}[stat]{Lemma}
\newtheorem{cor}[stat]{Corollary}
\newtheorem*{thm}{Theorem}
\theoremstyle{definition}

\newtheorem{mydef}[stat]{Definition}
\newcommand{\Pf}{\noindent{\bf Proof. }}

\newcommand{\hence}{\ensuremath{\Rightarrow}}
\newcommand{\isoto}{\xrightarrow{\sim}}

\newcommand{\DMO}{\DeclareMathOperator}
\DMO{\Pic}{Pic}
\DMO{\Cl}{Cl}
\DMO{\Proj}{Proj}
\DMO{\Spec}{Spec}
\newcommand{\SPec}{\mathop{\mathbf{Spec}}}

\renewcommand{\geq}{\geqslant}
\renewcommand{\leq}{\leqslant}

\newcommand{\LL}{\mathcal{L}}

\begin{document}
\maketitle
\begin{abstract}
Let $X_r$ be a smooth Del Pezzo surface obtained from $\PP^2$ by blow-up of
$r\leq 8$ points in general position. The multigraded group
$Cox(X_r)=\bigoplus_{L\in\Pic X_r} \G(X_r,L)$ can be almost canonically
endowed with a ring structure, called the Cox ring of $X_r.$ We show that 
in characteristic $0$ this ring is Gorenstein and its spectrum has rational
singularities. Also a gap in math.AG/0309111 is pointed out.
\end{abstract}

Let $X_r$ be a smooth Del Pezzo surface obtained from $\PP^2$ by blow-up of
$r\leq 8$ points in general position (\ie no line passes through any three
points, no conic through any six points and no cubic passing through eight
points has one of them for a double point). The integer $d=9-r$ is called
the the \emph{degree} of a Del Pezzo surface. It is known, cf.~\cite{Ma},
that the Picard group of $X_r$ equals $\ZZ^{r+1}$ and $d=(K,K)$ is the
self-intersection number of the canonical class of $X_r.$ The multigraded
group $Cox(X_r)=\bigoplus_{L\in\Pic X_r} \G(X_r,L)$ can be almost
canonically endowed with a ring structure, called \emph{the Cox ring} of
$X_r$ after similar constructions studied in \cite{Cox} and \cite{HK}.

This ring was studied in \cite{BaPo} together with its spectra
$\A(X_r):=\Spec Cox(X_r)$ and $\PP(X_r):=\Proj Cox(X_r),$ where in the
latter case the $\ZZ$-grading on the Cox ring is induced from the Picard
grading by taking intersection with the anticanonical class $-K$. It was
shown that for $r\geq 3$ this ring is generated by elements of degree 1,
namely, by global sections of the invertible sheaves whose divisors are
exceptional curves (one section for each curve) and, for $r=8$, also by two
global sections of the anticanonical sheaf. It was also shown that for
$r\leq 6$ this ring is defined by quadratic relations up to radical
(actually, the latter holds for $r\leq 7,$ but the last case was left
out of the paper in order to get a more elegant statement). These
properties showed the similarity of $Cox(X_r)$ to the homogeneous
coordinate ring of the orbit of the highest weight vector in the
projectivization of a certain representation of a certain algebraic group,
which can be naturally associated to a Del Pezzo surface, see \cite{BaPo}
for details.

In this paper we prove further properties of $Cox(X_r)$ that demonstrate
this similarity. After recalling the definition of the ring structure on
$Cox(X_r)$ and other useful facts about this ring and related objects, we
recall a description (originating in the work of Nagata) of the Cox ring as
a ring of invariants of a linear representation of a unipotent group, which
gives us for free that this ring is normal and factorial (it was shown in
\cite{BaPo} that it is finitely generated). Then we prove the main result:
\begin{thm}
In characteristic $0$ the variety $\A(X_r)$ has rational singularities and
is Gorenstein.
\end{thm}
Therefore in characteristic $0$ the Cox ring is Cohen-Macaulay. We also
calculate the anticanonical class of $\PP(X_r)$ (that is, the a-invariant
of the Cox ring) which equals $9-r$ times the ample generator of the Picard
group.

For homogeneous coordinate rings of orbits of the highest weight vector in
the projectivization of a representation of a certain algebraic group 
the rationality of singularities was shown in \cite{Ra}. As the Picard
group of the orbits of the highest weight vectors corresponding to Del
Pezzo surfaces is $\ZZ,$ their homogeneous coordinate rings are also
factorial and Gorenstein.

{\bf Remark on \cite{BaPo}.} Proposition~3.4 of \cite{BaPo}, which claims
that for $r=8$ the subring of $Cox(X_8)$ generated by global sections of
the invertible sheaves whose divisors are exceptional curves coincides with
$Cox(X_8)$ in every Picard-homogeneous component except the
$(-K)$-component, is supplied there with a sketch of a proof that can only
show that it is true for a generic $X_8.$ I therefore apologize and take
the responsibility of publishing a fallacious proof and state that this
claim is still open.

\medskip
I would like to thank V.~V.~Batyrev for this problem and the ideas used in
\cite{BaPo} and three more people for introducing me to the methods I use
here: I.~V.~Arzhantsev to invariant theory, D.~A.~Timashev to the
description of the Cox ring as a ring of invariants and D.~I.~Panyushev 
to rational singularities.

\section{Preliminaries}

We study the ring $Cox(X_r)=\bigoplus_{L\in\Pic(X_r)}\G(X_r,L)$ for a Del
Pezzo surface $X_r$ obtained as a blowup of $r$ points, $3\leq r\leq 8.$ We
are also interested in its spectrum $\A(X_r)$, its projective spectrum
$\PP(X_r)$ and in the sheafy variations of these concepts, as the next
lemma shows:

\begin{lemma}\label{sheaf}
Let $\mathcal{R}_{X_r}:=\bigoplus_{L\in\Pic(X_r)}L$ be a sheaf of
$\OO_{X_r}$-algebras and let $T$ be the torus with the character lattice
$\Pic(X_r).$ Then $\A(X_r)$ is the closure of a quasiaffine embedding of
the principal $T$-bundle $\myT(X_r):=\SPec\mathcal{R}_{X_r},$ where $\SPec$
denotes the relative spectrum.
\end{lemma}

{\bf Remark.} 
$T$ is the N\'eron-Severi torus \cite[5.3]{MaTs} of $X_r$ and $\myT(X_r)$
is a universal torsor in the sense of \cite[5.4]{MaTs}.
\medskip

Now we describe the multiplication in the sheaf; the proof of the lemma is
postponed.

The multiplication in the ring (and in the sheaf) is well-defined,
if we fix a basis in the Picard group and a particular representative $L_i,
i=1,\dots,k$ in each basis class: then in every class we choose a sheaf of
type $L^{\bar n}:=L_1^{\otimes n_1}\otimes\dots\otimes L_k^{\otimes n_k},
\ \bar n\in \ZZ^k,$ so that we have
natural bilinear maps of sections from two sheaves to their tensor product
and natural contractions $L^{\otimes n}\otimes L^{\otimes (-m)}
\to L^{\otimes(n-m)}$: as our sheaves have rank 1, the tensor
multiplication is commutative, so it doesn't matter which index is
contracted with which.

The sheaf, hence the ring, does not depend on the choice of this basis and
of representatives: let $M_1,\dots, M_k$ be
also representatives for a basis of the Picard lattice,
$\varphi_i\co M_i\isoto L^{\bar s_i}$
 isomorpisms of sheaves. Let $S$ be the matrix with rows $\bar s_i.$
Then let us take for each $M^{\bar m}$
the isomorphism $\varphi^{\bar m}\co M^{\bar m} \to L^{\bar mS},$ which is
the composition of $\bigotimes_i\varphi_i^{\otimes m_i}$ with the contractions
of sheaves $L_i$ with their inverses. Then we claim, that this family of
sheaf isomorphisms defines an isomorphism of the sheaf of algebras
constructed by means of $M_i$ to the sheaf of algebras constructed by means
of $L_j.$ The map is bijective and respects the addition, as we have sheaf
isomorphisms. What we need to show is that they respect the multiplication,
\ie the commutativity of the following diagram:

$$
\xymatrix{
M^{\bar m_1}\otimes M^{\bar m_2} \ar[r]\ar[d]^-{\varphi^{\bar 
m_1}\otimes\varphi^{\bar m_2}} 
& M^{\bar m_1+\bar m_2}\ar[d]^-{\varphi^{\bar m_1+\bar m_2}}\\ 
\left(\bigotimes_{i=1}^kL^{m_{1i}\bar s_i}\right)\otimes
\left(\bigotimes_{i=1}^kL^{m_{2i}\bar s_i}\right) \ar[r]\ar[d]&
\bigotimes_{i=1}^kL^{(m_{1i}+m_{2i})\bar s_i}.\ar[d]\\
L^{\bar m_1S}\otimes L^{\bar m_2S} \ar[r]& L^{(\bar m_1+\bar m_2)S}
}
$$

Here unlabeled arrows are suitable contractions. The lower square
commutes because the contractions commute with one another, and the upper
one does because the contractions are natural and commute with isomorphisms
of sheaves (where one takes the inverse of the adjoint isomorphism for the
dual sheaf). To put it short, the isomorphism evidently respects the
multiplication.

Another description of this ring is given in \cite[Def.~1.1]{BaPo}.

The ring $Cox(X_r)$ is graded by the semigroup of effective divisors, and
there is also a coarser grading given by
$Cox(X_r)_i=\bigoplus_{(L,-K)=i}\G(X_r,L),$ with respect to which we'll
speak of a degree of an element and of a divisor and consider the
projective spectrum $\PP(X_r).$ One sees that our ring is an integral
domain because all the homogeneous elements \wrt the fine grading are
nonzerodivisors.

{\bfseries Proof of Lemma~\ref{sheaf}.} The torus $T$ acts naturally on
$\myT(X_r),$ preserving the fibers, and, as each fiber is isomorphic to the
spectrum of the group algebra $\Bbbk[\Pic(X_r)],$ this is a principal
bundle. As $Cox(X_r)$ is the ring of global sections of $\OO_{\myT(X_r)}$ and
it is finitely generated, it determines a morphism of $\myT(X_r)$ into an
affine space, and the closure of the image is $\A(X_r).$ We need to show
that this morphism is an embedding. First of all, as the homogeneous
coordinate ring of the anticanonical embedding of $X_r$ is a subring of
$Cox(X_r),$ the morphism does not glue points in different fibers together
and its tangent map can degenerate only along the fibers. So we have to
check that it embeds each fiber. And for this it is enough to show that the
affine coordinate ring of the fiber is the localization of the image of
$Cox(X_r)$ in it. The map from $Cox(X_r)$ to the coordinate ring of the
fiber over a point $p,$ which is the group algebra of $\Pic(X_r),$ takes to
zero the global sections of $L\subset\mathcal R$ vanishing at $p$, and it
takes other sections to a scalar times the class of $L$ in $\Pic(X_r).$ We
need to show that the sheaves that have non-vanishing at $p$ global
sections generate the Picard group. But base-point free sheaves do generate
the Picard group by virtue of \cite[Chap.~II, Theorem~5.17]{Ha}.
\qed

Let us quote another useful lemma \cite[Prop.~4.4 and its proof]{BaPo}:
\begin{lemma}\label{cones}
Let $X_{r-1}$ be the Del Pezzo surface obtained by the contraction
$\pi\co X_r\to X_{r-1}$ of an exceptional curve $E$ on $X_r$, $x_E$ the
generator of $Cox(X_r)$ corresponding to $E$ and $U_E$ the affine chart
$\{x_E\ne 0\}$ on $\PP(X_r)$. Then there is a natural isomorphism
$$
U_E \cong \A(X_{r-1})
$$
which can be described as follows on the level of rings: it maps
the $n$-th homogeneous component $(Cox(X_{r-1}))_n$ to 
$\pi^*(Cox(X_{r-1}))_n/x_E^n\subset (Cox(X_r)_{x_E})_0.$ 
\end{lemma} 

We need one more result \cite[Theorem~3.2]{BaPo}: 
\begin{lemma}
For $3\leq r \leq 8$, the ring $Cox(X_r)$ is generated by elements of
degree 1. If $r\leq 7$, then the generators of $Cox(X_r)$ are global
sections of invertible sheaves defining the exceptional curves. If $r=8,$
then we should add to the above set of generators two linearly independent
global sections of the anticanonical sheaf on $X_8$.
\end{lemma}

\begin{cor}\label{charts} 
For $4\leq r\leq 8$ $\PP(X_r)$ is covered with affine charts isomorphic to
$\A(X_{r-1}),$ where the varieties $X_{r-1}$ are obtained from $X_r$ by
contraction of exceptional curves.
\end{cor}
\Pf
We conserve the notations of Lemma~\ref{cones}. For $4\leq r\leq 7$ the
previous lemma shows that the generators $x_E$ corresponding to the
exceptional curves generate the irrelevant ideal of the $\ZZ$-graded ring
$Cox(X_r).$ For $r=8$ we consider the subring $R$ of $Cox(X_8)$ generated
by $x_E$'s for all exceptional curves $E\subset X_8.$ This is a graded
subring and its projective spectrum $Z$ is covered with the charts
$V_E=\{x_E\ne 0\}\subset Z.$ The inclusion of graded rings induces a map
$$
j\co \bigcup_E U_E \to Z
$$ 
from an open subset of $\PP(X_8).$ We claim that this map is an
isomorphism. To show this it is enough to show that the maps
$j^{-1}(V_E)\to V_E$ are isomorphisms, and $j^{-1}(V_E)$ evidently equals
$U_E.$ The map $j$ induces isomorphisms $U_E\isoto V_E,$ because the
inclusion of graded rings induces the inclusion of affine coordinate rings
$\Bbbk[V_E]\subset\Bbbk[U_E]\cong Cox(X_7),$ but $Cox(X_7)$ is generated by
$x_{E'}$'s for $(E',E)=0,$ and the description of the isomorphism shows 
that $\Bbbk[U_E]$ is generated by $x_{E'}/x_E$'s, which lie in 
$\Bbbk[V_E],$ so this subring coincides with the whole ring. Now
an open set $\bigcup_E U_E$ of an irreducible projective variety 
$\PP(X_8)$ is itself isomorphic to a projective variety $Z.$ Therefore it 
coincides with $\PP(X_8).$

Thus in both cases the charts $U_E$ cover $X_r.$ And they are what we need
by Lemma~\ref{cones}.
\qed

\section{$Cox(X_r)$ as a Ring of Invariants}

D.~A.~Timashev made me aware of the following interpretation of $Cox(X_r),$
which is present in \cite[Sect.~3]{Na} for a blow-up of 16 points and is
analogous to the example in \cite{Na2}. I reproduce here, with his kind
permission, the statement and the proof from \cite[Lemmas 1 and 2]{Tim}.
\begin{stat}
Let a Del Pezzo surface $X_r$ be the blow-up of $r$ points
$P_j=(a_{0j}:a_{1j}:a_{2j})$ in $\PP^2.$ Then for the $2r$-dimensional
vector space $V$ with the action of the group
$$
U=\left\{
\begin{pmatrix}
\begin{matrix}
1 & u_1\\
0 & 1
\end{matrix}
& &\mbox{\LARGE $0$}\\
&\ddots&\\
\mbox{\LARGE $0$}& &
\begin{matrix}
1 & u_r\\
0 & 1
\end{matrix}
\end{pmatrix} \mbox{\LARGE$:$} \sum_{j=1}^ra_{ij}u_j=0,\ i=0,1,2
\right\}
$$
one has the ring of invariants $\Bbbk[V]^U\cong Cox(X_r).$
\end{stat}
\Pf 
Let us denote the coordinates on $V$ (\wrt the basis in which $U$ has the
above form) by $x_1,y_1,\dots,x_r,y_r.$ Let $\ang V$ be the open subset
$\{y_1\dots y_r\ne 0\},$ $z_j=x_j/y_j$ be coordinates on it, then
$\Bbbk[\ang V]=\Bbbk[y_j^{\pm 1},z_j]$ with $U$ acting by $y_j\mapsto y_j,$
$z_j\mapsto z_j+u_j.$ Then for $w_i=y_1\dots y_r\sum_ja_{ij}z_j$ one
obviously has $\Bbbk[\ang V]^U=\Bbbk[y_j^{\pm1},w_i].$ One also has
$\Bbbk[V]^U=\Bbbk[V]\cap\Bbbk[\ang V]^U.$ Both rings on the right admit a
grading by the total degree in $x$'s and a multigrading \wrt the $y$'s, so
we can consider multihomogeneous components $A_{d;m_1,\dots,m_r}$ of the
intersection, where $d$ stands for the degree in $w$'s and $m_j$ in
$y_j^{-1},$ namely, $A_{d;m_1,\dots,m_r}$ equals the set of
$F_d(w_0,w_1,w_2)y_1^{-m_1}\dots y_r^{-m_r}\in\Bbbk[V]$ for $F_d$
homogeneous of degree $d.$

The last condition means that $F_d$ as a polynomial in $x$'s and $y$'s is
divisible by $y_j^{m_j}.$ It is claimed that this in turn means that $F_d$
as a polynomial in $w$'s vanishes to order at least $m_j$ at (the blow-up
of) $P_j.$ We can check the claim for each $j$ separately, and to see this,
we notice that $GL_3$ acts on the equations and on the $w$'s, allowing us
to suppose that we check the claim at $P_1=(1:0:0),$ and that
$P_2=(0:1:0),$ $P_3=(0:0:1).$ Then $w_0\equiv x_1y_2\dots y_r
\mod{(x_4,\dots,x_r)},$ $w_1\equiv y_1x_2y_3\dots y_r
\mod{(x_4,\dots,x_r)},$ and $w_2\equiv y_1y_2x_3\dots y_r
\mod{(x_4,\dots,x_r)}.$ Also $w_1$ and $w_2$ are divisible by $y_1,$ and as
the vanishing order at $P_1$ is the maximal degree of the ideal
$(w_1,w_2)\in\Bbbk[w_1,w_2]$ that contains $F_d(1,w_1,w_2),$ $F_d(x,y)$ is
divisible by $y_1$ raised to this order. Conversely, if $F_d(x,y)$ is
divisible by $y_1^{m_1},$ then it remains divisible after substituting
$x_4=\dots=x_r=0,$ $x_1=y_2=\dots=y_r=1.$ The substitution takes $w_0$ to
$1,$ $w_1$ to $x_2y_1$ and $w_2$ to $x_3y_1,$ so the vanishing order is at
least $m_1.$

Thus we have identified the homogeneous components of $\Bbbk[V]^U$ with the
global sections $\G(\OO(dl_0-m_1l_1-\dots-m_rl_r))\subset \Bbbk(X_r),$ and
the multiplication in $\Bbbk[V]^U$ is componentwise inherited from
$\Bbbk(X_r).$ As the same description is valid for $Cox(X_r)$ (when the basis
$l_0,\dots,l_r$ is chosen in the Picard group), the rings are isomorphic.
\qed

\begin{cor}\label{lfac}
The ring $Cox(X_r)$ is normal and factorial for $3\leq r\leq 8.$ The
variety $\PP(X_r)$ is locally factorial and projectively normal and its
Picard group is $\ZZ$ generated by the hyperplane section from the
embedding for which $\A(X_r)$ is the affine cone.
\end{cor}
\Pf 
The facts about $Cox(X_r)$ are known (and easy) properties for the invariants
of a connected group without characters acting on a polynomial ring
\cite[Theorems 3.16 and 3.17]{AGIV}. Local factoriality follows from
the global one and other facts about $\PP(X_r)$ follow then from
\cite[Chap.~II, \S6, Exercise~6.3 (c)]{Ha}, which says that the affine cone
over a projective variety is factorial iff the variety is projectively
normal and has the divisor class group $\ZZ$ generated by the hyperplane
section.
\qed

\section{Rational Singularities}

In this section we assume that $\Bbbk$ is of characteristic zero and $r$ is
between 3 and 8. Under these restrictions we show the following proposition:

\begin{stat}\label{ratsing}
For $r\geq 3$ the ring $Cox(X_r)$ is Cohen-Macaulay and Gorenstein, the
singularities of $\PP(X_r)$ and $\A(X_r)$ are rational and the
anticanonical class of $\PP(X_r)$ is ample and equals $d$ times the ample
generator of the Picard group.
\end{stat}

{\bf Remark.} For $r\leq 2$ the Del Pezzo surface is a toric variety, so
its Cox ring is a polynomial ring by \cite{Cox} and still has rational 
singularities and is Gorenstein, so our theorem follows from the proposition.

\subsection{Preliminaries}

We recall from \cite{Ke} the notion of rational singularity and some of its
properties.

\begin{mydef}
An algebraic variety $X$ is said to have a \emph{rational singularity at}
$x\in X,$ if for some (and then for any) resolution of singularities
$p\co\T X\to X$ one has  $(Rp_*\OO_{\T X})_x=\OO_{X,x}.$ It is said to have
\emph{rational singularities,} if the previous condition is true at every
point.
\end{mydef}

This means that $p_*\OO_{\T X}=\OO_X,$ \ie $X$ is normal, and
$R^ip_*\OO_{\T X}=0$ for $i>0.$ In characteristic 0 the rational
singularity condition is equivalent to $X$ being Cohen-Macaulay and the
natural morphism $p_*\omega_{\T X} \to\omega_{X}$ being isomorphism for
some (and then any) resolution of singularities, where $\omega$ denotes the
canonical sheaf \cite{Ke}, by virtue of the following
\textbf{Grauert-Riemenschneider vanishing theorem} (cf.
\cite[Satz~2.3]{GR}):
\begin{lemma}
Let $p\co\T X\to X$ be a resolution of singularities (in characteristic 0),
$\omega_{\T X}$ the canonical sheaf of $\T X.$ Then $R^ip_*\omega_{\T X}=0$
for $i>0.$
\end{lemma}
Another useful fact is \textbf{Boutot's theorem} \cite[Th\'eor\`eme,
Cor.]{Bou}:
\begin{lemma}
Let $A\to B$ be a pure morphism of localisatons of affine rings over a field of
characteristic 0, \ie such that for each $A$-module $M$ the morphism $M\to
B\otimes_A M$ is injective. Then if the spectrum of $B$ has rational
singularities, the same is true for $A.$ In particular, the quotient of an
affine variety with rational singularities modulo the action of a reductive
group has rational singularities.
\end{lemma}

We cite a result which relates the rational singularity property of a
projective variety $X$ to that of its affine cone $\WH X.$

\begin{lemma}[an instance of {\cite[Theorem 1]{Ke2}}]\label{cone.rat.sing}
The affine cone $\WH X$ over a projective variety $X\subset \PP^n$ has rational
singularities iff $X$ is projectively normal, has rational singularities and
$H^i(X,\OO(k))=0$ for $i>0, k\geq 0.$ 
\end{lemma}

\Pf \cite[Theorem 1]{Ke2} states the following: let $X$ be a proper variety
with rational singularities, $L_1,\dots,L_n$ invertible sheaves on $X$
which are generated by global sections, $C=\Spec(\bigoplus_{m_j\geq
0}\G(X,\bigotimes L_j^{m_j}))$ be what \cite{Ke2} call a multicone. Then
$C$ has rational singularities iff $H^i(X,\bigotimes L^{m_j})=0$ for any
$i>0$ and any $m_1\geq 0,\dots, m_n\geq 0$ and $H^i(X,\bigotimes
L^{-1-m_j})=0$ for any $i<\dim X$ and any $m_1\geq 0,\dots, m_n\geq 0.$ In
characteristic 0 we can drop the last condition as it is shown in the proof
of \cite[Theorem~1]{Ke2} to be equivalent to the Grauert--Riemenschneider
vanishing theorem.

For the ``\hence" direction we notice that $X$ has rational singularities
by Boutot's theorem, being a quotient of $\WH X\setminus\{0\}$ by the
$\Bbbk^*$-action. Also the cone must be normal by definition of the
rational singularity, so $X$ is projectively normal and $\WH
X=\Spec(\bigoplus_{n\geq 0}\OO_X(n)),$ so we can apply
\cite[Theorem~1]{Ke2} and obtain our vanishing condition.

For the inverse direction we say that $X$ is projectively normal, so
$$
\WH X=\Spec\left(\bigoplus_{n\geq 0}\OO_X(n)\right)
$$ 
and we can apply 
\cite[Theorem~1]{Ke2}. \qed

The following lemma on vanishing will be useful in checking the conditions of
the previous lemma:

\begin{lemma}[an instance of the main result of {\cite{GR}}]\label{Kodaira}
The Kodaira vanishing (\ie $D$ ample $\hence H^i(D+K)=0$ for $i>0$) holds on a
projective variety with rational factorial Gorenstein singularities.
\end{lemma}

\Pf The above mentioned main result \cite[Einleitung, see also Satz~2.1]{GR}
states that $H^i(X, K_X\otimes V)=0$ for $i>0,$ where $X$ is an irreducible
complex Moishezon space (also a projective variety will do), $K_X=p_*(K_{\T
X})$ for a resolution of singularities $p\co\T X \to X$ is the Grauert
canonical sheaf, so in case of rational singularities it coincides with the
Grothendieck canonical sheaf $\omega_X$ which is invertible because the
singularities are Gorenstein, and $V$ is a positive vector bundle, so an
ample invertible sheaf will do.
\qed

\subsection{The Anticanonical Class of $\PP(X_r)$.}

First we deal with some geometric properties of $\PP(X_r).$ 

\begin{lemma}
Let $T'$ be the torus with the character lattice $\langle -K\rangle^{\perp}
\subset \Pic(X_r)$. Then the natural principal $T'$-bundle $\myT'(X_r)$ is an
open subset of $\PP(X_r)$ and this embedding induces an isomorphism of the
Picard groups.
\end{lemma}

\Pf In the notations of Lemma~\ref{sheaf} let us consider the projection
$\myT(X_r)\to \PP(X_r).$ As the projection from the affine cone without apex to
the corresponding projective variety is open, the image is an open subset.

This projection is the factorization (principal bundle) of $\myT(X_r)$
modulo the action of a 1-dimensional subtorus $\GG_m$ of the N\'eron-Severi
torus $T,$ which has the property that all the characters of $T$
corresponding to exceptional curves restrict to the identity character of
$\GG_m.$ Therefore the map of character groups $\WH T\to \WH{\GG_m}$ may be
described as $\Pic(X_r)\ni \chi\mapsto (\chi,-K)\in \ZZ.$ The duality
between tori and their character groups gives us then the claimed
description of the quotient torus $T'.$ As the map $\myT(X_r)\to X_r$ is
constant on $T$-orbits, it is constant on $\GG_m$-orbits, so it factorizes
through $\myT'(X_r),$ and the map $\myT'(X_r)\to X_r$ is still the
factorization modulo the $T$-action, which pushes down to a $T'$-action and
thus gives the structure of a principal $T'$-bundle.

There is also a principal $\GG_m$-bundle $\myT(X_r)\to\myT'(X_r).$ For a
principal bundle one has \cite[Chap.~I, \S3]{Mum}
$\Pic(\myT'(X_r))\cong\Pic^{\GG_m}\myT(X_r),$  and, as the only invertible
global regular functions on $\myT(X_r)$ are $Cox(X_r)^*=\Bbbk^*,$ different
characters of $\GG_m$ define different linearizations of the structure
sheaf on $\myT(X_r)$ and one has an exact sequence
$$
0\to\WH{\GG_m}\to\Pic^{\GG_m}\myT(X_r)\to\Pic\myT(X_r)=0,
$$ 
where the last equality follows from the fact that $\myT(X_r)$ is an open
subset of the spectrum of a factorial ring $Cox(X_r).$ Thus one sees that
the Picard group of $\myT'(X_r)$ is $\ZZ,$ as well as its divisor class
group, because this variety is smooth. By virtue of Cor.~\ref{lfac} these
groups for $\PP(X_r)$ also equal $\ZZ,$ as $\PP(X_r)$ is locally factorial.
Now the statement about the Picard (or the divisor class) groups follows,
because by~\cite[Chap.~II, \S6, Prop.~6.5 a)]{Ha} we have a surjection
$$
\Cl(\PP(X_r))\to\Cl(\myT'(X_r))\to 0
$$ 
with both groups isomorphic to $\ZZ,$ so it is an isomorphism.\qed

Now we calculate the anticanonical class of $\myT'(X_r).$
\begin{lemma}\label{antican}
The anticanonical class of $\myT'(X_r)$ equals $9-r$ times the ample
generator of $\Pic(\myT'(X_r)).$
\end{lemma}

\Pf Let $\pi\co\myT'(X_r)\to X_r$ be the natural projection, $TM$ be the
tangent bundle of a variety $M.$ Then, as with every principal bundle,
there is an exact sequence of vector bundles on $\myT'(X_r)$:
$$
0\to\LL(T')\times\myT'(X_r)\to T\myT'(X_r)\to\pi^*TX_r\to 0,
$$  
where the inclusion map is the partial derivative of the action map \wrt the
group element.

As both varieties $X_r$ and $\myT'(X_r)$ are smooth, we have Chern classes of
vector bundles on them, so from the triviality of the first bundle in the
sequence above one sees that the Chern classes of (the tangent bundle of)
$\myT'(X_r)$ are the pullbacks of the Chern classes of $X_r.$ In particular, it
is true for the anticanonical class $c_1.$ One sees from the discussion
above that the map $\Pic(X_r)\to \Pic \myT'(X_r)\cong \ZZ$ in the intersection
with $-K,$ so $-K$ is mapped to $d=9-r$ times a generator, and for sure the
ample one, as pullback respects the effective cone.\qed

We'll see later, that the anticanonical class of $\PP(X_r)$ is well-defined and
coincides with that of $\myT'(X_r),$ whence the title of this section.

\subsection{Proof of the Proposition}  
We proceed by induction on $r,$ the case $r=3$ being evident from
$\PP(X_3)=\PP^5$ \cite[Prop.~3.1]{BaPo}. The case $r=4,$ with
$\PP(X_4)=G(3,5)$ \cite[Prop.~4.1]{BaPo}, the Grassmann variety of 3-planes
in a 5-dimensional vector space, is also known (see, e.g., \cite{Ra}), but
we would not call it evident.

We introduce two series of statements:
\begin{itemize}
\item[$(A_r)$] 
The ring $Cox(X_r)$ is Gorenstein and a rational singularity.
\item[$(B_r)$] 
The variety $\PP(X_r)$ has rational Gorenstein singularities
and an anticanonical class as claimed by the Proposition.
\end{itemize}

Now we do a two-step induction.

``$(A_r)\hence(B_{r+1})$":  Corollary~\ref{charts} states that for $r\leq
7$ the variety $\PP(X_{r+1})$ is covered with charts isomorphic to
$\A(X_r)$. Therefore it has rational Gorenstein singularities, in
particular, its (anti)canonical sheaf is invertible. As the inclusion map
of $\myT'(X_{r+1})$ evidently pulls back the anticanonical sheaf of
$\PP(X_{r+1})$ to the one of $\myT'(X_{r+1}),$ the claim follows from
Lemma~\ref{antican}.

``$(B_r)\hence(A_{r})$": first we show that $Cox(X_r)$ is a rational
singularity. By virtue of Lemma~\ref{cone.rat.sing} one needs only to check
the vanishing of the higher cohomologies of $\OO_{\PP(X_r)}(n)$ for
nonnegative $n,$ and this follows from Lemma~\ref{Kodaira} and the
ampleness of the anticanonical class.

Now, as rational singularities are Cohen-Macaulay \cite{Ke} and factorial
Cohen-Macaulay rings which have a canonical module (e.g., finitely
generated algebras) are Gorenstein \cite[Exercise~21.21]{Ei}, $Cox(X_r)$ is
Gorenstein. \qed


\begin{thebibliography}{Ke-Ra}
\bibitem[B-P]{BaPo}
V.~V.~Batyrev, O.~N.~Popov. The Cox Ring of a Del Pezzo Surface. \pp To
appear in: Arithmetic of higher-dimensional algebraic varieties, B.~Poonen
and Y.~Tschinkel (eds.), Progress in Math. 226 (2004), Birkh\"auser. Also
math.AG/0309111.
%
\bibitem[Bou]{Bou}
\begin{otherlanguage}{french}
J.-F.~Boutot. Singularit\'es rationnelles et quotients par les groupes
r\'eductifs. \pp Invent. Math. 88, 65--68 (1987). 
\end{otherlanguage}
%
\bibitem[Cox]{Cox}
D.~Cox. The Homogeneous Coordinate Ring of a Toric Variety.
\pp J. Algebr. Geom. 4, No.1, 17--50 (1995).
%
\bibitem[Ei]{Ei}
D.~Eisenbud. Commutative Algebra with a view towards Algebraic Geometry.
Graduate Texts in Mathematics, 150. Springer-Verlag, 1995.
%
\bibitem[G-R]{GR}
\begin{otherlanguage}{german}
H.~Grauert, O.~Riemenschneider. Verschwindungss"atze f"ur analytische
Kohomologiegruppen auf komplexen R"aumen  \pp Invent. Math. 11, 263--292
(1970).
\end{otherlanguage}
%
\bibitem[H-K]{HK} Y.~Hu, S.~Keel. Mori Dream Spaces and GIT. \pp Michigan 
Math. J. {\bf 48}, 331--348 (2000). Also math.AG/0004017.
%
\bibitem[Ha]{Ha} R.~Hartshorne. Algebraic Geometry.
Graduate Texts in Mathematics, 52. Springer-Verlag, 1977.
%
\bibitem[Ke]{Ke}
G. Kempf. Cohomology and Convexity. \pp G.~Kempf, F.~Knudsen, D.~Mumford,
B.~Saint-Donat. Toroidal embeddings. I. Springer Lect. Notes Math. 339. Chapter
I, \S3, 49--52 (1973).
%
\bibitem[Ke-Ra]{Ke2}
G. Kempf,  A. Ramanathan. Multi-cones over Schubert Varieties. \pp Invent.
Math. 87, 353--363 (1987).
%
\bibitem[Ma]{Ma}Yu.~I.~Manin.
Cubic Forms. Algebra, Geometry, Arithmetic. North-Holland Mathematical Library.
Vol. 4. North-Holland, 1974.
%
\bibitem[Ma-Ts]{MaTs}
Yu.~I.~Manin, M.~A.~Tsfasman. Rational varieties: Algebra, geometry and
arithmetic. \pp Russ. Math. Surv. 41, No.2, 51--116 (1986). Translation
from Usp. Mat. Nauk 41, No.2(248), 43--94 (1986).
%
\bibitem[Mum]{Mum}
D.~Mumford. Geometric invariant theory. 
\begin{otherlanguage}{german}
Ergebnisse der Mathematik und ihrer Grenzgebiete. Neue Folge. 34. 
Springer-Verlag, 1965. 
\end{otherlanguage}
%
\bibitem[Na]{Na}
M.~Nagata. On the fourteenth problem of Hilbert \pp Proceedings of the 
international congress of mathematicians (14--21 August 1958). Cambridge, 
1960. P.~459--462.
%
\bibitem[Na2]{Na2}
M.~Nagata. On the 14th problem of Hilbert \pp Amer. J. Math. 81, 766--772 
(1959).
%
\bibitem[Ra]{Ra} 
A. Ramanathan. Schubert varieties are arithmetically Cohen-Macaulay \pp
Invent. Math. 80, No. 2, 283--294 (1985).
%
\bibitem[Ti]{Tim}
D.~A.~Timashev. The 14th problem of Hilbert and rational surfaces (in
Russian),

{\ttfamily
http://www.math.msu.su/department/algebra/staff/timashev/rat-surf.ps }
%
\bibitem[V-P]{AGIV}
E.~B.~Vinberg, V.~L.~Popov. Invariant theory. \pp A.~N.~Parshin,
I.~R.~Shafarevich, R.~V.~Gamkrelidze (eds.). Algebraic geometry IV.
Encyclopaedia of mathematical sciences. 55. Berlin: Springer-Verlag, 1994.
%
\end{thebibliography}
\end{document}